\documentclass[11pt,a4paper]{article}
\usepackage{epsfig}
\usepackage{amssymb}
\usepackage{amsthm}
\usepackage{amsmath}
\usepackage{latexsym}
\usepackage{epic}

\newcommand{\mc}{\mathcal}
\newcommand{\pro}{\mc P}
\newcommand{\nb}{\mathbf{n}}
\newcommand{\kb}{\mathbf{k}}
\newcommand{\ab}{\mathbf{a}}
\newcommand{\ph}{{\tilde{p}}}
\newcommand{\ep}{\epsilon}
\newcommand{\bep}{b_{\epsilon}(d)}
\newcommand{\Ex}{\mathrm{E}}
\newcommand{\Var}{\mathrm{Var}}
\newcommand{\Cov}{\mathrm{Cov}}
\newcommand{\beqrn}{\begin{eqnarray}}
\newcommand{\eeqrn}{\end{eqnarray}}
\newcommand{\beq}{\begin{equation}}
\newcommand{\eeq}{\end{equation}}
\newcommand{\qt}{\tilde{q}}
\newcommand{\Qt}{\tilde{Q}}
\newcommand{\weak}{\stackrel{\rm D}{\to}}
\newcommand{\prob}{\stackrel{\rm P}{\to}}
\newcommand{\nbinf}{\mbox{ as } \nb \to \infty}
\newcommand{\Or}{{\rm O}}
\newcommand{\Op}{{\rm O}_{\rm P}}

\newcommand{\oor}{{\rm o}}

\newcommand{\lsq}{q}
\newcommand{\lsqep}{\tilde{q}_{\ep}}
\newcommand{\Lsq}{{\tilde{Q}}_{\nb}}

\theoremstyle{plain}
\newtheorem{theorem}{Theorem}
\newtheorem{lemma}{Lemma}

\newtheorem{proposition}[theorem]{Proposition}

\newenvironment{reflist}{\begin{list}{}
         {\itemsep=20pt \parsep=3pt
          \topsep=0pt  \parskip=20pt  \listparindent=-.15in
           \leftmargin= 0.15in }
         \item \ \vspace{-.35in} }
         {\end{list}}

\begin{document}
\author{David K\"{a}llberg\thanks{Corresponding author. 
 E-mail address:  {\tt david.kallberg@math.umu.se} }, \; Oleg\ Seleznjev \vspace{0.5cm}\hfill \\
Department of Mathematics and Mathematical Statistics\\
Ume{\aa } University, SE-901 87 Ume\aa , Sweden\\
}
\date{}
\title{ Estimation of entropy-type integral functionals }
\maketitle

\begin{abstract}
Entropy-type integral functionals of densities are widely used in mathematical statistics, 
information theory,
and computer science.
Examples include measures of closeness between distributions (e.g., density power divergence) 
and uncertainty characteristics for a random variable (e.g., R\'{e}nyi entropy). 
In this paper, we study $U$-statistic estimators for a class of such functionals.
The estimators are based on $\ep$-close vector observations in the corresponding independent and identically distributed samples. 
We prove asymptotic properties of the estimators 
(consistency and asymptotic normality) 
under mild integrability and smoothness conditions for the densities.
The results can be applied in diverse problems in mathematical statistics and computer science
(e.g., distribution identification problems, approximate matching for random databases, two-sample problems).
\end{abstract}
 \baselineskip=3.4 ex
\noindent \emph{AMS 2010 subject classification:} 62G05, 62G10, 62G20, 94A17
\medskip


\section{Introduction}
Let the distributions $\mc P_X$ and $\mc P_Y$
of the $d$-dimensional random variables $X$ and $Y$
have densities $p_X(x)$ and $p_Y(x), x \in R^d$, respectively.
Various characteristics in mathematical statistics, information theory,
and computer science, say \textit{entropy-type integral functionals}, 
are expressed in terms of integrated polynomials in $p_X(x)$ and $p_Y(x)$.
For example, a widely accepted measure of closeness between
$\mc P_X$ and $\mc P_Y$ is the (quadratic) density power divergence (Basu et al., 1998)
\beq
D_2=D_2(\mc P_X,\mc P_Y) := \int_{R^d} (p_X(x)-p_Y(x))^2 dx. \notag
\eeq
Other examples include the (integer order) R\'{e}nyi entropy for quantifying uncertainty in $X$ (R\'{e}nyi, 1970)
\beq
h_k=h_k(\mc P_X) := \frac{1}{1-k}\log \left( \int_{R^d}p_X(x)^k dx \right), \quad k=2,3,\ldots, \notag
\eeq
and the differential variability for some database problems (Seleznjev and Thalheim, 2010)
\beq
v=v(\mc P_X, \mc P_Y) := -\log \left( \int_{R^d}p_X(x)p_Y(x) dx \right). \notag
\eeq
Henceforth we use $\log x$ to denote the natural logarithm of $x$.
For non-negative integers $k_1,k_2 \geq 0,\kb :=(k_1,k_2)$,
we consider the R\'{e}nyi entropy functionals (K\"{a}llberg et al., 2012)
\beq
q_\kb=q_{k_1,k_2} :=  \int_{R^d} p_X(x)^{k_1} p_Y(x)^{k_2} dx, \quad k_1+k_2 \geq 2.  \notag
\eeq
Moreover, given a set of constants $\ab:= \{a_0,a_1,a_2\}$,
we introduce the related quadratic functionals
\beq
\lsq = \lsq(\ab) := a_{0} q_{2,0} + a_{1} q_{1,1} + a_{2} q_{0,2}.\notag
\eeq
Note that the quadratic divergence $D_2 = q_{2,0} -  2q_{1,1}+ q_{0,2}$,
the R\'{e}nyi entropy  $h_k = \log(q_{k,0})/(1-k), k=2,3,\ldots$, and the variability $v = -\log(q_{1,1})$.
Some applications of R\'{e}nyi entropy and divergence measures
can be found, e.g., in information theoretic learning (Principe, 2010).
More applications of entropy and divergence in statistics
(e.g., distribution identification problems and statistical inference),
computer science (e.g., average case analysis for random databases,
pattern recognition, and image matching), and econometrics are discussed,
e.g., in Pardo (2006), Broniatowski et al.\ (2012), Leonenko et al.\ (2008), Escolano et al.\ (2009),
Seleznjev and Thalheim (2003, 2010), Thalheim (2000), Leonenko and Seleznjev (2010), Neemuchwala et al.\ (2005), and Ullah (1996).
The divergence $D_2$ belongs to a subclass of the
Bregman divergences that find various applications in statistics
(see, e.g., Basseville, 2010, and references therein).

In this paper, to demonstrate the general approach,
we study estimation of some entropy-type integral functionals, e.g., $q_\kb$ and $\lsq(\ab)$, using independent samples from $\mc P_X$ and $\mc P_Y$.
We prove asymptotic properties for a class of $U$-statistic estimators for these functionals.
The estimators are based on $\ep$-close vector observations in the corresponding samples.
We generalize some results and techniques contained in Leonenko and Seleznjev (2010)
 and K\"{a}llberg et al.\ (2012).
In particular, we obtain consistency of the estimators under more general density conditions and prove asymptotic normality of the estimators for the quadratic functionals $\lsq({\ab})$.

Leonenko et al.\ (2008) study asymptotic properties of nearest-neighbor estimators for $q_\kb$
and establish consistency when the densities are bounded. 
Ahmad and Cerrito (1993) obtain asymptotic normality of a kernel estimate of the quadratic divergence $D_2$, but only under quite strong differentiability conditions for the densities. 
In the one-dimensional case, Bickel and Ritov (1988) and Gin\'{e} and Nickl (2008) show rate optimality, efficiency, and asymptotic normality of kernel-type estimators for $q_{2,0}$.
Laurent (1996) uses orthogonal projection to build an efficient and asymptotically normal estimator of $q_{2,0}$
for multivariate distributions.

The number of pairs of $\ep$-close observations (or the number of small inter-point distances) in a random sample is among the most studied examples of $U$-statistics with kernels varying with the sample size (see, e.g., Weber, 1983, Jammalamadaka and Janson, 1986, Penrose, 1995). 
A significant feature of this and related statistics is their asymptotically normality under weak assumptions.
In our work, we extend this property to the two-sample case and construct asymptotically normal estimators of the quadratic functionals $q(\ab)$ under mild smoothness conditions for the densities.  
In particular, the proposed method is valid for densities of low regularity (or smoothness) given that a rate of convergence slower than $\sqrt{n}$ is acceptable.
This appears to be a clear advantage of our approach compared to similar studies, see, e.g., Bickel and Ritov (1988), Gin\'{e} and Nickl (2008).

First we introduce some notation.
Let $d(x,y):=|x-y|$ denote the Euclidean distance
in $R^d$ and define $B_\ep(x):=\{y:d(x,y)<\ep\}$ to be an open
$\ep$-ball in $R^d$ with center at $x$ and radius $\ep$.
Denote by $\bep:=\ep^db_1(d), b_1(d)=2\pi^{d/2}/(d\Gamma(d/2))$, the volume of the $\ep$-ball.
Define the $\ep$-ball probability as
\beq
p_{X,\ep}(x):=P\{X \in B_\ep(x)\}. \notag
\eeq
We say that the vectors $x$ and $y$ are $\ep$\textit{-close} if $d(x,y) < \ep$ for some $\ep >0$.
Let $X_1,\ldots,X_{n_1}$ and $Y_1,\ldots,Y_{n_2}$
be mutually independent samples of independent and identically distributed
(i.i.d.) observations from $\pro_X$ and $\pro_Y$, respectively.
Define $\nb:= (n_1,n_2)$, $n:=n_1+n_2$, and say that $\nb \to \infty$ if  $n_1, n_2 \to \infty$.
Throughout the paper, we assume that $\ep = \ep(\nb) \to 0  \nbinf$.

Denote by $\weak$ and $\prob$ convergence in distribution and probability, respectively.
For a sequence of random variables $\{V_\nb\}$, we write $V_\nb = \Op(1) \nbinf$ if for any $\delta > 0$
and large enough $n_1,n_2 \geq 1$, there exists $C>0$ such that $P(|V_\nb|>C) \leq \delta$.
Moreover, for a numerical sequence $\{u_\nb\}$, let $V_\nb = \Op(u_\nb) \nbinf$ if $V_\nb/u_\nb = \Op(1) \nbinf$.

In what follows, we consider estimation problems for both one and two samples.
However, in the statements of results and the proofs, it is assumed, for sake of space and clarity,
that two samples are available, i.e.\ $n_1,n_2 \geq 0$.
This can be done without loss of generality, because in the one-sample case,
e.g., estimation of $q_{k_1,0}, k_1 \geq 2$, from $X_1,\ldots, X_{n_1}$,
an auxiliary sample $Y_1, \ldots, Y_{n_2}$ can be used.

The main goal of this paper is to provide asymptotic properties for estimation under weak distributional assumptions. We leave for further research important matters such as efficiency of the estimators under study, optimality of the obtained convergence rates, and selection of an optimal parameter $\ep=\ep(\nb)$.

The remaining part of the paper is organized as follows.
In Section \ref{sec:qk}, we consider estimation of the R\'{e}nyi entropy functional $q_\kb$.
Section \ref{sec:q2}  presents the asymptotic properties for estimation of the quadratic functional $\lsq(\ab)$. In Section \ref{sec:app}, we consider applications of the results to estimation of divergence,
the two-sample problem, and statistical inference for some entropy-type characteristics.
Numerical examples illustrate the rate of convergence of the asymptotic results.
Section \ref{sec:proofs} contains the proofs of the statements from the previous sections.

\section{Estimation of the R\'{e}nyi entropy functional $q_{\kb}$}\label{sec:qk}
We introduce the $U$-statistic estimators of $q_\kb$ proposed by K\"{a}llberg et al.\ (2012).
For non-negative integers $r$ and $m$, define $\mc S_{m,r}$ to be the set of all $r$-subsets of $\{1,\ldots, m\}$.
Let $S \in \mc S_{n_1,k_1}$ and $T \in \mc S_{n_2,k_2}$.
When $k_1 \geq 1$, we define
\beq
\psi^{(i)}_{\kb,\nb,\ep}(S;T) =I(d(X_i,X_j)<\ep, d(X_i,Y_l) <\ep, \forall j \in S, \forall l \in T), \quad i \in S,  \notag \\
\eeq
i.e.\ the indicator of the event that the observations $\{X_j, j \in S\}$ and $\{ Y_l, l\in T\}$ are $\ep$-close to $X_i$. In a similar way, for $k_1 = 0$ and $k_2 \geq 1$, let
\beq
\psi^{(i)}_{\kb,\nb,\ep}(T) =I(d(Y_i,Y_j) <\ep, \forall j \in T), \quad i \in T. \notag \\
\eeq
By conditioning, we have
\beq
q_{\kb,\ep} :=
\left \{
\begin{array}{ll}
\Ex ( \psi^{(i)}_{\kb,\nb}(S;T)) = \Ex (p_{X,\ep}(X)^{k_1-1} p_{Y,\ep}(X)^{k_2}), & \mbox{ if } k_1 \geq 1, \\
\Ex ( \psi^{(i)}_{\kb,\nb}(T)) = \Ex (p_{Y,\ep}(Y)^{k_2-1} ), & \mbox{ if } k_1 = 0, k_2 \geq 1.
 \end{array}
 \right.
 \notag
\eeq
Now, a $U$-statistic for $q_{\kb,\ep}$ (see, e.g., Ch.\ 2, Lee, 1990) is given by
\beq
Q_{\kb,\nb} = Q_{\kb,\nb,\ep} : = \binom{n_1}{k_1}^{-1} \binom{n_2}{k_2}^{-1} \sum_{S \in \mc S_{n_1,k_1}} \sum_{T \in \mc S_{n_2,k_2}} \psi_{\kb,\nb}(S;T), \notag
\eeq
with the \textit{kernel} $\psi_{\kb,\nb}(S;T)$ defined by the symmetrization
\beq
 \psi_{\kb,\nb}(S;T) = \psi_{\kb,\nb,\ep}(S;T) :=
\left \{
\begin{array}{ll}
\displaystyle \frac{1}{k_1} \sum_{i \in S} \psi^{(i)}_{\kb,\nb,\ep}(S;T), & \mbox{ if } k_1 \geq 1, \\
\displaystyle \frac{1}{k_2} \sum_{i \in T} \psi^{(i)}_{\kb,\nb, \ep}(T), & \mbox{ if } k_1 = 0, k_2 \geq 1.
 \end{array}
 \right.
 \notag
\eeq
By definition, $Q_{\kb,\nb}$ is an unbiased estimator of $q_{\kb,\ep}$. Let $k:=k_1+k_2, k\geq 2$, and define the estimator of $q_\kb$ as
\beq \Qt_{\kb,\nb} := Q_{\kb,\nb}/\bep^{k-1}. \notag
\eeq
The asymptotic properties of $\Qt_{\kb,\nb}$ depend on the rate of decrease of $\ep(\nb)$.
In our first theorem, we establish consistency under two different asymptotic rates of $\ep(\nb)$ 
with minimal density assumptions in {\it (ii)}. 
Note that K\"{a}llberg et al.\ (2012) prove {\it (i)} under stronger density conditions 
(i.e.\ boundedness and continuity).  
\begin{theorem} Let ${n_1/n \to \rho, 0<\rho<1}$.
\label{th:main1}
\begin{itemize}
\item[(i)]
If $p_X,p_Y \in L_{2k-1}(R^d)$ and $n\ep^{d(1-1/k)} \to \infty$, then
$$
\Qt_{\kb,\nb} \prob q_\kb \nbinf.
$$
\item[(ii)]
If $p_X,p_Y \in L_{k}(R^d)$ and $n\ep^{d(k-1)} \to \infty$, then
$$
\Qt_{\kb,\nb} \prob q_\kb \nbinf.
$$
\end{itemize}
\end{theorem}

\section{Estimation of the quadratic functional $\lsq(\ab)$}\label{sec:q2}
The following linear combination is considered as an estimator of the quadratic functional $\lsq$,
\beq
\Lsq=\Lsq(\ab) := a_0 \Qt_{2,0,\nb} + a_1 \Qt_{1,1,\nb} + a_2 \Qt_{0,2,\nb}. \notag
\eeq
Theorem \ref{th:main1} gives conditions for consistency of the estimator $\Lsq$.
Next we describe some asymptotic normality properties of $\Lsq$. 
Let $\qt_{\kb,\ep} := \Ex(\Qt_{\kb,\nb}) = \bep^{-1} q_{\kb,\ep}$ and define
\beq
\lsqep := \Ex(\Lsq) = a_0 \qt_{2,0,\ep} + a_1 \qt_{1,1,\ep} + a_2 \qt_{0,2,\ep}. \notag
\eeq
We also introduce the characteristics
\begin{align}
\zeta = \zeta(\ab,\rho, \mc P_X, \mc P_Y) & := \frac{4}{\rho} \Var \left(a_{0}p_X(X) + \frac{a_{1}}{2}p_Y(X) \right) + \frac{4}{1-\rho} \Var \left(a_{2}p_Y(Y) + \frac{a_{1}}{2}p_X(Y) \right), \notag \\
\eta = \eta(\ab,\rho, \mc P_X, \mc P_Y) & :=  \frac{2}{ b_1(d)}\bigg(\frac{a_{0}^2}{\rho^2} q_{2,0} + \frac{a_{2}^2}{(1-\rho)^2} q_{0,2} + \frac{a_{1}^2}{2 \rho (1 - \rho)} q_{1,1}\bigg ), \quad 0<\rho < 1. \notag
\end{align}
Henceforth, in order to have $\eta > 0$, we exclude the (trivial) cases $a_0=a_2=a_1 q_{1,1}=0$.
\begin{theorem}\label{th:main2} Let $p_X,p_Y \in L_3(R^d)$ and $n_1/n = \rho, 0< \rho < 1$.
\begin{itemize}
\item[(i)] If $n\ep^d \to \beta, 0 < \beta\leq\infty,$ (and $\zeta > 0$ when $\beta = \infty$), then
\beq
\sqrt{n}(\Lsq - \lsqep) \weak N(0, \zeta + \eta/\beta) \nbinf. \notag
\eeq
\item[(ii)] If $n\ep^d \to 0$ and $n^2\ep^d \to \infty$, then
\beq
n\ep^{d/2}(\Lsq - \lsqep) \weak N(0,\eta) \nbinf. \notag
\eeq
\end{itemize}
\end{theorem}
To ensure a sufficient rate of decay for the bias of the estimator $\Lsq$, we propose smoothness (or regularity) conditions for the densities.
Denote by $H_{2}^{(\alpha)}(K), 0<\alpha \leq 1, K>0$, a linear space of functions in $L_3(R^d)$
that satisfy an $\alpha$-H\"{o}lder condition in $L_2$-norm with constant $K$, i.e.\ if $p \in H_2^{(\alpha)}(K)$ and $h \in B_1(d)$, then
\beq\label{L2}
||p(\cdot+h)- p(\cdot)||_2 \leq K|h|^{\alpha}.
\eeq
Note that \eqref{L2} holds if, e.g., for some function $g \in L_2(R^d)$,
\beq
|p(x+h)-p(x)| \leq g(x) |h|^{\alpha} \notag
\eeq
and hence $H^{(\alpha)}_2(K)$ is wider than, e.g., the corresponding H\"{o}lder class considered by Bickel and Ritov (1988).
There are different ways to introduce the density smoothness, e.g.,
by the pointwise H\"{o}lder conditions (K\"{a}llberg et al., 2012) or the Fourier characterization (Gin\'{e} and Nickl, 2008).

In the next theorem, we present a bound for the bias and the rate of convergence in probability of $\Lsq$ in terms of the density smoothness $\alpha$.
 Let $L(n), n \geq 1$, be a slowly varying function as $n \to \infty$.
\begin{theorem}\label{th:main3}  Let $p_X,p_Y \in H^{(\alpha)}_2(K)$ and $n_1/n \to \rho, 0<\rho<1$.
\begin{itemize}
\item[(i)] Then for the bias, we have
\beq
|\lsqep-\lsq|\leq C \ep^{2\alpha}, C>0. \notag
\eeq
\item[(ii)] If $0<\alpha \leq d/4$ and $\ep \sim cn^{-1/(2\alpha + d/2)}, c > 0$, then
\beq
\Lsq - \lsq = \Op (n^{ -2\alpha/(2\alpha+d/2)}) \nbinf. \notag
\eeq
\item[(iii)] If $\alpha > d/4$, $\ep \sim L(n)n^{-1/d}$, and $n\ep^d \to \beta, 0<\beta \leq \infty$, then
\beq
\Lsq - \lsq = \Op (n^{ -1/2}) \nbinf. \notag
\eeq
\end{itemize}
\end{theorem}
\medskip
To make the asymptotic normality results of Theorem~\ref{th:main2} practical
(e.g., to construct approximate confidence intervals), the unknown asymptotic variances have to be estimated. For this, we need consistent estimators of $\zeta$ and $\eta$.
By expanding the terms in $\zeta$, we see that it is a
 function of $\rho$ and the functionals $\{q_{i,j}: 2 \leq i+j  \leq 3\}$,
i.e.\ $\zeta = \zeta(\rho,\{q_{i,j}: 2 \leq i+j  \leq 3\})$.
If $p_X,p_Y \in L_3(R^d)$ and the sequence $\ep_0 = \ep_0(\nb)$ satisfies $n\ep_0^{2d} \to \infty$, then Theorem \ref{th:main1}{\it (ii)} yields that $\{\Qt_{i,j,\nb,\ep_0}: 2 \leq i+j  \leq 3 \}$
are consistent estimators of these functionals.
Hence, we set up a plug-in estimator of $\zeta$ according to
\beq
\zeta_\nb=\zeta_{\nb,\ep_0} := \zeta(\rho_\nb,\{\Qt_{i,j,\nb,\ep_0}: 2 \leq i+j  \leq 3 \}), \notag
\eeq
where $\rho_\nb:=n_1/n$.
Denote by $v^2_\nb = v^2_{\nb,\ep_0}$ the corresponding estimator for $\eta$ and let
$w^2_{\nb} = w^2_{\nb,\ep_0} :=  \zeta_\nb + v^2_\nb/(n\ep^d)$
be an estimate of $\zeta + \eta/ \beta$ when $n\ep^d \to \beta, 0<\beta \leq \infty$.

Now we combine the results of Theorem \ref{th:main2}, the smoothness conditions, and variance estimators to construct asymptotically pivotal quantities. These can be used,  e.g., for calculating asymptotic confidence intervals for the functional $\lsq$.
It is worth noting that normal limits are possible under any density smoothness $0<\alpha \leq 1$.
However,
observe that the obtained rate of convergence is slower than $\sqrt{n}$ in the 'low regularity case' $\alpha \leq d/4$.
\begin{theorem}\label{th:main4}  Let $p_X,p_Y \in H^{(\alpha)}_2(K)$ and $n_1/n = \rho, 0<\rho<1$.
\begin{itemize}
\item[(i)] If $\alpha > d/4$, $\ep \sim L(n)n^{-1/d}$, and $n\ep^d \to \beta, 0<\beta \leq \infty$,
(and $\zeta > 0$ when $\beta = \infty$), then
\beq
\sqrt{n}(\Lsq - \lsq) \weak N(0, \zeta + \eta/\beta) \quad   \mbox{and} \quad
\sqrt{n}(\Lsq - \lsq)/w_\nb \weak N(0, 1) \nbinf. \notag
\eeq
\item[(ii)] If $\alpha > (d/4) \gamma$, for some $0<\gamma<1$, and $\ep \sim c n^{-2/((1+\gamma)d)}, c>0$, then
\beq
n^{\gamma/(1+\gamma)}c^{d/2}(\Lsq - \lsq) \weak N(0, \eta) \quad   \mbox{and}  \quad
n^{\gamma/(1+\gamma)} c^{d/2}(\Lsq - \lsq)/v_\nb \weak N(0, 1) \nbinf. \notag
\eeq
\item[(iii)] If $\ep \sim L(n)^{2/d}n^{-2/d}$ and $L(n) \to \infty$, i.e.\ $n^2\ep^d \to \infty$, then
\beq
L(n)(\Lsq - \lsq) \weak N(0, \eta) \quad \mbox{and} \quad
L(n)(\Lsq - \lsq)/v_\nb \weak N(0, 1) \nbinf. \notag
\eeq
\end{itemize}
\end{theorem}
\noindent
The generality of this result can be formulated in the following way:
given a lower bound $\alpha^*$ for the smoothness $\alpha$, we get from {\it (i)} and {\it (ii)} that an asymptotically normal estimator of $q$ is available with $n^{\nu}$ rate of convergence,
where $0<\nu \leq 1/2$ is a function of $\alpha^*$ and the dimension $d$.
Furthermore, with no information about $\alpha$, ${\it (iii)}$ implies that asymptotically normal estimation of $q$ is still possible, but at the slower rate $L(n)$.
\medskip

\noindent \textbf{Remark 1.} {\it (i)}
The condition $n_1/n=\rho, 0 < \rho < 1,$ in Theorems \ref{th:main2} and \ref{th:main4} is technical
and we claim that it can be replaced with the slightly weaker condition $n_1/n \to \rho, 0 < \rho < 1$.
\medskip

\noindent
{\it (ii)}
Notice that the condition $\zeta > 0$ in Theorems 2 and 4 is unnecessary if we let $n\ep^d \to \beta, 0\leq \beta < \infty$.
In general, this condition imposes some restrictions on the densities.
For example, when considering the divergence $D_2$,  $\zeta>0$ implies that $p_X(x) \neq p_Y(x)$ on a set of positive measure. For the functional $q_{2,0}$, we have $\zeta > 0$ unless $X$ is uniformly distributed on some set $\mc D \subset R^d$ (Jammalamadaka and Janson, 1986).
\medskip

\noindent
{\it (iii)}
In the one-sample case, the results of Theorem \ref{th:main4} (and Theorem 2) are essentially independent of $\rho$.
In fact, consider, e.g., the estimator $\Qt_{2,0,\nb}$ of $q_{2,0}$, i.e.\ $\ab = \{1,0,0 \}$, $\zeta = 4\Var(p_X(X))/\rho$,
and $\eta = 2b_1(d)^{-1}q_{2,0}/\rho^2$.
We have $n = n_1/\rho$, so if $n_1 \ep^d \to \lambda, 0<\lambda < \infty$,
then $n\ep^d \to \lambda/\rho =: \beta$.
Hence, it follows from Theorem \ref{th:main4}\textit{(i)} that
\beq
\sqrt{n_1}(\Qt_{2,0,\nb}-q_{2,0}) \weak
N\left(0, 4\Var(p_X(X)) + \frac{2}{b_1(d)}q_{2,0}/\lambda \right) \mbox{ as } n_1 \to \infty. \notag
\eeq
Therefore, we obtain a result with $\sqrt{n_1}$-scaling that does not depend on $\rho$ as desired.
Similar modifications can be done for the scalings in {\it (ii)} and {\it (iii)}.

\section{Applications}\label{sec:app}
\subsection{Estimation of divergence}
The introduced quadratic divergence $D_2$ belongs to
the wide class of density power divergences (Basu et al., 1998), defined by
\begin{align}
D_s &=D_s(\mc P_X, \mc P_Y) := \frac{1}{s-1}q_{s,0} - \frac{s}{s-1}q_{1,s-1} + q_{0,s}, \quad s > 1. \notag
\end{align}
Another family of divergences is the pseudodistances introduced in Broniatowski et al.\ (2012)
\begin{align}
R_s  =
R_s(\mc P_X, \mc P_Y)
 := \frac{1}{s}\log \left(q_{s,0} \right)
- \frac{1}{s-1}\log \left(q_{s-1,1} \right)+ \frac{1}{s(s-1)}\log \left(q_{0,s} \right)
, \quad s>1. \notag
\end{align}
Note that, for non-negative integers $r = 2,3,\ldots$,
Theorem~\ref{th:main1} enables the construction of consistent plug-in estimators
$\hat{D}_{r,\nb}$ and $\hat{R}_{r,\nb}$ of the divergences $D_r$ and $R_r$, respectively.
Moreover, the quadratic estimator $\hat{D}_{2,\nb}$ is asymptotically normal under the conditions of Theorem~\ref{th:main4} (see also Remark 1{\it (ii)}).

The quadratic divergence $D_2$ can be used as a dissimilarity measure to investigate pairwise differences among $M$ populations or objects.
Let the features of population $l$ be represented by the random feature vector $V_l$
with density $p_{V_l}(x) , x \in R^d, l=1,\ldots, M$.
Using independent samples from populations $ V_l, l=1,\ldots, M $, we apply, e.g.,
the Bonferroni method in combination with Theorem~\ref{th:main4} and
 calculate the $\binom{M}{2}$ approximate simultaneous confidence intervals $\{ I_{lm} \}$
for the quadratic divergences $\{ D_{2,l,m} \}$, $D_{2,l,m} := D_2(\mc P_{V_l},\mc P_{V_m})$, for a given confidence level.
These intervals can be used to determine which populations are different with respect to their feature densities.
\bigskip

\begin{figure}[htbp]
\begin{center}
\includegraphics[width=0.85\textwidth,height=0.35\textheight]{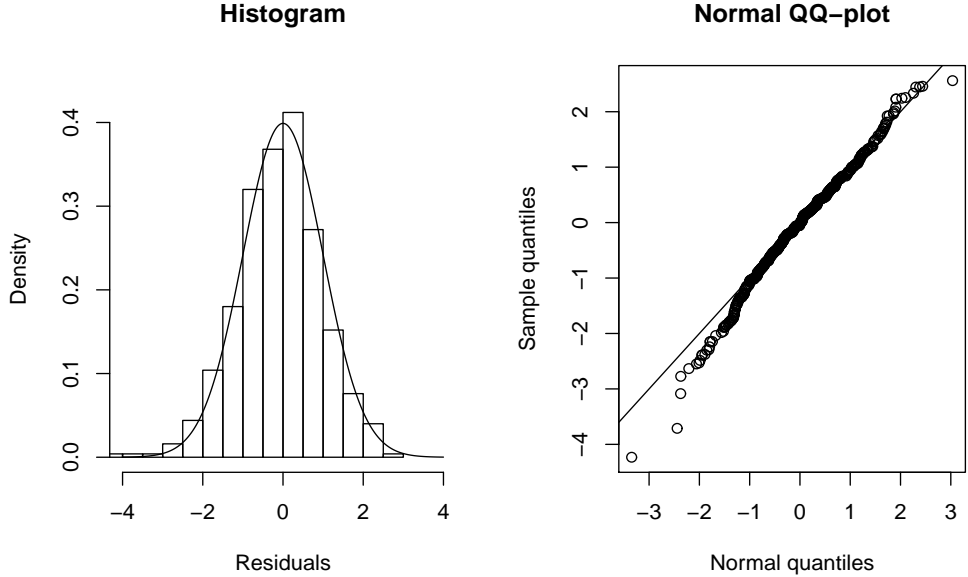}
\caption{
Three-dimensional distributions with $t(3)$-i.i.d.\ and $N(1,1)$-i.i.d.\ components, respectively; sample sizes $n_1=n_2=500$
and $\ep = \ep_0= 1/4$. Standard normal approximation for the normalized residuals; $N_{sim} = 500$.}
\label{fig1}
\end{center}
\end{figure}

\noindent
\textbf{Example 1.} We consider estimation of the quadratic density power divergence $D_2(\mc P_X,\mc P_Y)$ between two three-dimensional distributions.
The distribution of the components of $X$ and $Y$ are $t(3)$-i.i.d.\ and $N(1,1)$-i.i.d., respectively.
In this case it holds that $D_2 \approx 0.018$.
We simulate $N_{sim} = 500$ pairs of independent samples from $\mc P_X$ and $\mc P_Y$ and calculate the corresponding normalized residuals $R^{(i)}_\nb:=\sqrt{n}(\hat{D}_{2,\nb}-D_2)/w_\nb$, $i =1,\ldots,N_{sim}$,
with $n_1=n_2=500$, and $\ep=\ep_0 =1/4$.
The histogram and normal quantile plot in Figure~\ref{fig1} illustrate the performance of the normal approximation of $R^{(i)}_\nb$ indicated by Theorem~\ref{th:main4}{\it (i)}.
The p-value (0.41) for the Kolmogorov-Smirnov test also fails to reject the hypothesis of standard normality of the residuals.
\subsection{The two-sample problem}
A general null hypothesis of closeness between $\mc P_X$ and $\mc P_Y$ is given by
\beq
H_0: p_X(x) = p_Y(x) \mbox{ a.e. } \notag
\eeq
We consider the problem of testing $H_0$ against the alternative
$H_1$ that $p_X(x)$ and $p_Y(x)$ differ on a set of positive measure (often referred to as the two-sample problem).
Note that the alternative can be written as $H_1: D_2 > 0$. Hence, we define a test statistic based on the estimator $\hat{D}_{2,\nb}$ for $D_2$ (see, e.g., Li, 1996) according to
\beq
T_\nb:= \frac{n\ep^{d/2}} {v_\nb} \hat{D}_{2,\nb}. \notag
\eeq
The next proposition presents the asymptotics for the distribution of $T_\nb$. 
Li (1996) proves a similar result (for a general kernel) under more restrictive density conditions 
(i.e.\ boundedness and continuity) (see also Ahmad and Cerrito, 1993). Let $\{c_\nb\}$ be a numerical sequence such that $c_\nb = \oor(n\ep^{d/2}) \nbinf$.
\begin{proposition}\label{test} Assume that $p_X,p_Y \in L_3(R^d)$, $n^2\ep^d \to \infty$, and $n_1/n = \rho, 0<\rho<1$.
\begin{itemize}
\item[(i)] Under $H_0$, we have
\beq
n\ep^{d/2}\hat{D}_{2,\nb} \weak N(0,\eta) \quad  \mbox{and} \quad T_\nb \weak N(0,1) \nbinf. \notag
\eeq
\item[(ii)] Under $H_1$, we have
\beq
P(T_\nb > c_\nb) \to 1 \nbinf. \notag
\eeq
\end{itemize}
\end{proposition}
\noindent
Thus, we reject $H_0$ if $T_\nb> \lambda_a$, where $\lambda_a$ is the $a-$quantile of the standard normal distribution.
This test has asymptotic significance level $a$
and is consistent against all alternatives that satisfy $p_X,p_Y \in L_3(R^d)$.
\medskip

\noindent
\textbf{Remark 2.} If $n\ep^d \to \infty$, we claim that the condition $p_X,p_Y \in L_3(R^d)$ in
Proposition~\ref{test} can be weakened to the (minimal) assumption $p_X,p_Y \in L_2(R^d)$
provided that the sequence $\ep_0=\ep_0(\nb)$ in $v^2_\nb = v^2_{\nb,\ep_0}$ satisfies $n\ep_0^d \to \infty $ (e.g., $\ep_0$ = $\ep$).
\subsection{Estimation of R\'{e}nyi entropy and differential variability}
Consider the class of functionals
\beq
h_\kb=h_{\kb}(\mc P_X, \mc P_Y) := \frac{1}{1-k}\log(q_\kb)
, \quad k \geq 2. \notag
\eeq
For $\mc P_X = \mc P_Y$, we have the R\'{e}nyi entropy $h_{k,0}$, which is one of
a family of functions for measuring uncertainty or randomness of a system (R\'{e}nyi, 1970).
Another important example is the differential variability $v=h_{1,1}$ 
for some database problems (Seleznjev and Thalheim, 2010).
When the densities are bounded and continuous, the results in K\"{a}llberg et al.\ (2012)
imply consistency of the truncated plug-in estimator
\beq
H_{\kb,\nb} := \frac{1}{1-k} \log (\max(\Qt_{\kb,\nb},1/n)). \notag
\eeq
for $h_\kb$. It follows from Theorem \ref{th:main1} that $H_{\kb,\nb}$
is consistent under weaker (integrability) conditions for the densities.

For the quadratic case $k=2$, i.e.\ $\kb = (2,0), (1,1)$, the asymptotic normality properties of $H_{\kb,\nb}$ are studied by
Leonenko and Seleznjev (2010) and K\"{a}llberg et al.\ (2012). The following proposition generalizes some of these results (see also Remark 1).

\begin{proposition}\label{log}
Assume that $k=2$. Let $p_X,p_Y \in H_{2}^{(\alpha)}(K)$ and $n_1/n=\rho$, $0<\rho<1$.
\begin{itemize}
\item[(i)]
If $\alpha > d/4$, $\ep \sim L(n)n^{-1/d}$, and $n\ep^d \to \beta, 0<\beta \leq \infty$,
(and $\zeta > 0$ when $\beta = \infty$), then
\beq
\sqrt{n}\,\Qt_{\kb,\nb}(H_{\kb,\nb} - h_\kb)/w_\nb \weak N(0, 1) \nbinf. \notag
\eeq
\item[(ii)] If $\alpha > (d/4) \gamma$, for some $0<\gamma<1$, and $\ep \sim c n^{-2/((1+\gamma)d)}, c>0$, then
\beq
n^{\gamma/(1+\gamma)}c^{d/2}\Lsq(H_{\kb,\nb} - h_\kb)/v_\nb \weak N(0, 1) \nbinf. \notag
\eeq
\item[(iii)] If $\ep \sim L(n)^{2/d}n^{-2/d}$ and $L(n) \to \infty$, i.e.\ $n^2\ep^d \to \infty$, then
\beq
L(n)\Qt_{\kb,\nb}(H_{\kb,\nb} - h_\kb)/v_\nb \weak N(0, 1) \nbinf. \notag
\eeq
\end{itemize}
\end{proposition}
\medskip
\noindent
The estimator $H_{\kb,\nb}$ of $h_\kb$ can be used, e.g., for distribution identification problems and approximate matching in stochastic databases
(for a description, see  K\"{a}llberg et al.,\ 2012).
\bigskip

\begin{figure}[htbp]
\begin{center}
\includegraphics[width=0.85\textwidth,height=0.35\textheight]{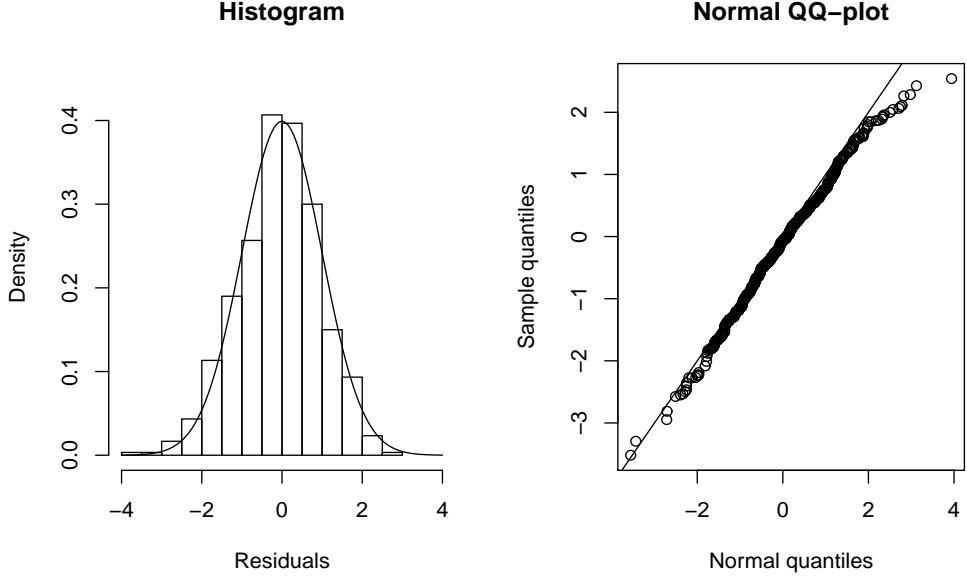}
\caption{Uniform distributions, $U(0,1)$ and $U(0,\sqrt{2})$; sample sizes $n_1=n_2=300$ and $\ep = \ep_0 = 1/100$.
Standard normal approximation for the normalized residuals; $N_{sim} = 600$.
}
\label{fig2}
\end{center}
\end{figure}

\noindent
\textbf{Example 2.}
Let $X$ and $Y$ be one-dimensional uniform random variables, i.e.\ $X \sim U(0,1)$ and $Y \sim~U(0,\sqrt{2})$,
and consider estimation of the differential variability $v=h_{1,1} = \log(2)/2$.
We simulate samples from $\mc P_X$ and $\mc P_Y$ and obtain the residuals $R^{(i)}_\nb:=\sqrt{n}\Qt_{1,1,\nb}(H_{1,1,\nb}-h_{1,1})/w_\nb, i=1,\ldots, N_{sim}$,
with $n_1=n_2=300, \ep=\ep_0 = 1/100$, and $N_{sim} = 600$.
Figure \ref{fig2} illustrates the normal approximation for these residuals indicated by Proposition \ref{log}\textit{(i)}.
The histogram, quantile plot, and p-value (0.36) for the Kolmogorov-Smirnov test
support the hypothesis of standard normality.

\section{Proofs}
\label{sec:proofs}
The following lemma is used in the subsequent proofs.
\begin{lemma}\label{lemma}
For $a,b \geq 0$, assume that $p_X, p_Y \in L_{a+b+1}(R^d)$. Then
\begin{equation*}
\bep^{-(a+b)}\Ex (p_{X,\ep}(X)^a p_{Y,\ep}(X)^b) \to q_{a+1,b}  \mbox{ as } \ep \to 0.
\end{equation*}
\end{lemma}
\noindent
{\it Proof.} Let $\ph_{X,\ep}(x) := \bep^{-1}p_{X,\ep}(x), \ph_{Y,\ep}(x) := \bep^{-1}p_{Y,\ep}(x)$.
Consider the decomposition
\begin{align}\label{lemma1:dc}
\bep^{-(a+b)}\Ex(p_{X,\ep}(X)^a &  p_{Y,\ep}(X)^b)   =  \int_{R^d} \ph_{X,\ep}(x)^a \ph_{Y,\ep}(x)^b p_X(x) dx \notag \\
 =&\int_{R^d} p_X(x)^{a+1}p_Y(x)^bdx + \int_{R^d} (\ph_{Y,\ep}(x)^b-p_Y(x)^b)p_X(x)^{a+1}dx \\
&+ \int_{R^d} ( \ph_{X,\ep}(x)^a-p_X(x)^a) \ph_{Y,\ep}(x)^b p_X(x) dx. \notag
\end{align}
We see that the assertion follows if the last two terms in \eqref{lemma1:dc}  tend to 0 as $\ep \to 0$.
By the extension of H\"{o}šlder's inequality (see, e.g., Ch.\ 2, Bogachev, 2007), we have
\begin{align} \label{lemma1:holder}
\Big |\int_{R^d}  (\ph_{X,\ep}(x)^a  -  p_X(x)^a) &  \ph_{Y,\ep}(x)^b  p_X(x) dx\Big |  \\
             & \leq ||\ph_{X,\ep}(\cdot)^a-p_X(\cdot)^a||_{(a+b+1)/a} ||\ph_{Y,\ep}(\cdot)^b||_{(a+b+1)/b} ||p_X(\cdot)||_{a+b+1}. \notag
\end{align}
The Lebesgue differentiation theorem implies
\beq\label{lemma1:diff}
\ph_{X,\ep}(x)^{a+b+1} \to p_X(x)^{a+b+1} \mbox{ as } \ep \to 0 \mbox{ a.e. }
\eeq
If $V = (V_1,\ldots, V_d)'$ is an auxiliary random vector uniformly distributed in the unit ball
$B_1(d)$, then $\ph_{X,\ep}(x) = \Ex ( p_X(x-\ep V))$ and thus Jensen's inequality leads to
\begin{align}\label{lemma1:dom1}
(\ph_{X,\ep}(x)^a)^{(a+b+1)/a} \leq g_\ep (x):=\Ex (p_X(x-\ep V)^{a+b+1})  = \frac{1}{\bep} \int_{B_\ep(x)}p_X(y)^{a+b+1}dy.
\end{align}
Since $p_X \in L_{a+b+1}(R^d)$, the Lebesgue differentiation theorem gives
\beq\label{lemma1:dom2}
g_\ep(x) \to g(x):= p_X(x)^{a+b+1} \mbox{ as } \ep \to 0 \mbox{ a.e. }
\eeq
Furthermore, Fubini's theorem yields
\beq\label{lemma1:dom3}
\int_{R^d} g_\ep(x)dx = \int_{R^d} g(x) dx.
\eeq
From \eqref{lemma1:diff}-\eqref{lemma1:dom3}
and a generalization of the dominated convergence theorem (see, e.g., Ch.\ 2, Bogachev, 2007), we get that
\beq\label{lemma1:conv1}
||\ph_{X,\ep}(\cdot)^a ||_{(a+b+1)/a} \to ||p_{X}(\cdot)^a ||_{(a+b+1)/a} \mbox{ as } \ep \to 0.
\eeq
In a similar way, we obtain
\beq\label{lemma1:conv2}
||\ph_{Y,\ep}(\cdot)^b||_{(a+b+1)/b} \to ||p_{Y}(\cdot)^b||_{(a+b+1)/b}\mbox{ as }\ep \to 0.
\eeq
Now we use the following result (see, e.g., Ch.\ 1, Kallenberg, 1997):
for a sequence of functions $f_n \in L_p(R^d), p \geq 1,n=1,\ldots$, with $f_n(x) \to f(x)$ a.e., $f \in L_p(R^d)$, it holds that
\beq\label{lemma1:Lp}
||f_n||_p \to ||f||_p \mbox{ iff } ||f_n - f||_p \to 0 \mbox{ as } n \to \infty.
\eeq
Note that \eqref{lemma1:diff}, \eqref{lemma1:conv1}, and \eqref{lemma1:Lp} imply
\beq\label{lemma1:conv3}
||\ph_{X,\ep}(\cdot)^a - p_X(\cdot)^a||_{(a+b+1)/a} \to 0 \mbox{ as } \ep \to 0.
\eeq
Finally, it follows from \eqref{lemma1:holder}, \eqref{lemma1:conv2}, and \eqref{lemma1:conv3} that
\beq
\int_{R^d} ( \ph_{X,\ep}(x)^a-p_X(x)^a) \ph_{Y,\ep}(x)^b p_X(x) dx \to 0 \mbox{ as } \ep \to 0. \notag
\eeq
By a similar argument, this also holds for the second term in \eqref{lemma1:dc}.
This completes the proof.
\hfill $\Box$ \\\\\\
\noindent \textit{Proof of Theorem \ref{th:main1}.}
For $l = 0,\ldots, k_1,$ and $m = 0, \ldots, k_2$, let
\begin{align}
  \psi_{\kb,l,m,\nb}(x_1,  \ldots, & x_l;y_1,\ldots,y_m) \notag  \\
                        & :=  \Ex(\psi_{\kb,\nb}(x_1,\ldots,x_l,X_{l+1},\ldots,X_{k_1};y_1,\ldots,y_m,Y_{m+1},\ldots,Y_{k_2}))  \notag
\end{align}
and
\beq
\sigma^2_{\kb,l,m,\nb} :=\Var(\psi_{\kb,l,m,\nb}(X_1,\ldots,X_l;Y_1,\ldots,Y_m)), \notag
\eeq
where we define $\sigma^2_{\kb,0,0,\nb}=0$.
From the conventional theory of $U$-statistics (see, e.g., Ch.\ 2, Lee, 1990), we have
\begin{equation}  \label{eq:var}
\Var(\Qt_{\kb,\nb}) =
\bep^{-2(k-1)}\sum_{l=0}^{k_1}\sum_{m=0}^{k_2}\frac{\binom{k_1}{l}\binom{k_2}{m}\binom{n_1-k_1}{k_1-l}\binom{n_2-k_2}{k_2-m}}
{\binom{n_1}{k_1}\binom{n_2}{k_2}} \sigma^2_{\kb,l,m,\nb}.
\end{equation}
Without loss of generality, we assume that $k_1 \geq 1$.
Following the argument in K\"{a}llberg et al.\ (2012), it is straightforward to show
\beq\label{s2lm}
\sigma^2_{\kb,l,m,\nb} \leq \Ex ( p_{X,3\ep}(X)^{2k_1-l-1}p_{Y,3\ep}(X)^{2k_2-m})
\eeq
and so Lemma \ref{lemma} yields
\beq \label{eq:bound}
\sigma^2_{\kb,l,m,\nb} = \mbox{O}(\bep^{2k-l-m-1}) \nbinf.
\eeq
Using the condition $n_1/n \to \rho, 0<\rho < 1$, we obtain
\begin{align}  \label{eq:varlim}
\bep^{-2(k-1)}
\frac{\binom{k_1}{l}\binom{k_2}{m}\binom{n_1-k_1}{k_1-l}\binom{n_2-k_2}{k_2-m}}
{\binom{n_1}{k_1}\binom{n_2}{k_2}} &  \sigma^2_{\kb,l,m,\nb}  \sim C_{l,m} \frac{\bep^{-(2k-l-m-1)}\sigma^2_{\kb,l,m,\nb}} {n^{l+m}\ep^{d(l+m-1)}} \nbinf
\end{align}
for some constant $C_{l,m}> 0$.
For $r=1, \ldots, k$, we have
\beq\label{ned}
n^{r}\ep^{d(r-1)} = (n \ep^{d(1-1/r)})^{r} \geq(n \ep^{d(1-1/k)})^{r}.
\eeq
Now it follows from \eqref{eq:var}-\eqref{ned} and the condition $n\ep^{d(1-1/k)} \to \infty$ that
\beq\label{eq:conv}
\Var(\Qt_{\kb,\nb}) = \mbox{O}\left(\frac{1}{n\ep^{d(1-1/k)}} \right) \to 0 \nbinf. 
\eeq
Finally, from Lemma \ref{lemma} we get $\Ex (\Qt_{\kb,\nb}) = \bep^{-(k-1)}q_{\kb,\ep} \to q_\kb$,
which together with \eqref{eq:conv} implies the statement. \\\\
\textit{(ii)} The argument is similar to that of {\it (i)}, so we show the main steps only.
Note that bound \eqref{s2lm} and Lemma 1 give
\beq
\sigma^2_{\kb,l,m,\nb} \leq\Ex ( p_{X,3\ep}(X)^{k_1-1}p_{Y,3\ep}(X)^{k_2}) = \mbox{O}(\bep^{k-1}) \nbinf \notag
\eeq
and hence the condition $n\ep^{d(k-1)} \to \infty$ yields
\beq
\bep^{-2(k-1)}
\frac{\binom{k_1}{l}\binom{k_2}{m}\binom{n_1-k_1}{k_1-l}\binom{n_2-k_2}{k_2-m}}
{\binom{n_1}{k_1}\binom{n_2}{k_2}}   \sigma^2_{\kb,l,m,\nb} \sim C_{l,m} \frac{\bep^{-(k-1)}\sigma^2_{\kb,l,m,\nb}} {n^{l+m}\ep^{d(k-1)}} \to 0 \nbinf \notag
\eeq
for some $C_{l,m}>0$. By combining this with \eqref{eq:var}, we obtain $\Var(\Qt_{\kb,\nb}) \to 0$ and, since also in this case $\Ex (\Qt_{\kb,\nb}) = \bep^{-(k-1)}q_{\kb,\ep} \to q_\kb$, the assertion follows.  This completes the proof.
\hfill $\Box$
\\\\\\
\noindent
{\it Proof of Theorem \ref{th:main2}.}
If  the case $\beta=\infty$ is excluded, then {\it (i)} and {\it (ii)} can be expressed together as follows:
if $n\ep^d \to \beta, 0 \leq \beta < \infty$ and $n^2\ep^d \to \infty$, then
\beq\label{nep^d/2}
n\ep^{d/2}(\Lsq- \lsqep) \weak N(0, \eta + \beta\zeta) \nbinf.
\eeq
We only present the proof for \eqref{nep^d/2} since the argument is similar for the remaining case $\beta = \infty$ in {\it (i)}.
If $n_3 = n_3(\nb)$ is defined to be the greatest common divisor of $n_1$ and $n_2$,
then $n_1 = n_3 l$ and $n_2 = n_3 m$, where $l$ and $m$ are positive integers that satisfy $l/(l+m) = \rho$.
Consider the following \emph{pooled} random vectors in $R^{d(l+m)}$
$$
Z_i := (X_{l(i-1)+1},\ldots, X_{l i}, Y_{m(i-1)+1},\ldots, Y_{mi}), \quad i = 1,\ldots, n_3.
$$
The method of proof relies on the decomposition
\beq\label{jj:dec}
n \ep^{d/2}(\Lsq- \lsqep) =  n\ep^{d/2}\binom{n_3}{2}^{-1}\bep^{-1} (U_{\nb} - \Ex (U_{\nb})) + R_\nb,
\eeq
where $U_{\nb}$, to be defined below, essentially is a \textit{one-sample} $U$-statistic with respect to the i.i.d.\ sample $\{Z_1, \ldots, Z_{n_3} \}$.
The idea is to prove that the remainder $R_{\nb}$ tends to 0 in probability
and use the corresponding result from Jammalamadaka and Janson (1986) to show asymptotic normality for the first term in \eqref{jj:dec}.

For $z_i:=(x_{l(i-1)+1},\ldots, x_{l i},y _{m(i-1)+1},\ldots, y_{mi}), i = 1, \ldots, n_3$, introduce the kernels
\begin{align}
\phi_\nb^{(1)}(z_i,z_j) & :=  \sum_{s=1}^{l}\sum_{t=1}^{l} I\left(d\left(x_{l(i-1)+s},x_{l(j-1)+t}\right)<\ep \right), \notag\\
 \phi_\nb^{(2)}(z_i,z_j) & :=  \sum_{s=1}^{m}\sum_{t=1}^{m} I\left(d\left(y_{m(i-1)+s},y_{m(j-1)+t}\right)<\ep \right), \notag \\
\phi_\nb^{(3)}(z_i,z_j) & :=  \sum_{s=1}^{l}\sum_{t=1}^{m} I \left(d \left(x_{l(i-1)+s},y_{m(j-1)+t} \right)<\ep \right) \notag
+ \sum_{s=1}^{l}\sum_{t=1}^{m} I \left(d \left(x_{l(j-1)+s},y_{m(i-1)+t} \right)<\ep \right). \notag
\end{align}
Moreover, define
\begin{align}\label{jj:fngn}
f_{\nb}(z_i,z_j)  :=& \   a_{0}l^{-2} \phi_\nb^{(1)} (z_i,z_j) + a_{2}m^{-2} \phi_\nb^{(2)} (z_i,z_j) + a_{1}(2lm)^{-1} \phi_\nb^{(3)}(z_i,z_j), \notag \\
\mu_\nb :=& \  \Ex(f_\nb(Z_1,Z_2)) = \bep \lsqep,  \notag \\
g_{\nb} (z_i)     :=& \    \Ex ( f_{\nb}(z_i,Z_j) ) - \mu_\nb \\
              =& \   \frac{a_{0}}{l} \sum_{s=1}^{l} p_{X,\ep}(x_{l(i-1)+s}) + \frac{a_{2}}{m} \sum_{s=1}^{m} p_{Y,\ep}(y_{m(i-1)+s}) \notag \\
             &+ \frac{a_{1}}{2}\left( \frac{1}{l}\sum_{s=1}^{l} p_{Y,\ep}(x_{l(i-1)+s}) +\frac{1}{m}\sum_{s=1}^{m} p_{X,\ep}(y_{m(i-1)+s}) \right) - \mu_\nb.  \notag
\end{align}
Let
\begin{align}
M_\nb  := \sum_{i<j} I(d(X_i,X_j) <\ep), \quad V_\nb := \sum_{i<j} I(d(Y_i,Y_j) <\ep), \quad
W_\nb  := \sum_{i,j} I(d(X_i,Y_j) <\ep), \notag
\notag
\end{align}
and note that
\beq\label{jj:QnA}
\bep \Lsq = a_{0}\binom{n_1}{2}^{-1} M_\nb + a_{2}\binom{n_2}{2}^{-1} V_\nb + a_{1}(n_1 n_2)^{-1} W_\nb.
\eeq
Now consider  the decompositions
$$
M_\nb = M_\nb^{(1)} + M_\nb^{(2)}, \quad V_\nb = V_\nb^{(1)} + V_\nb^{(2)}, \quad W_\nb = W_\nb^{(1)} + W_\nb^{(2)},
$$
where
\begin{align}
M_\nb^{(1)}  :=\sum_{i<j} \phi_\nb^{(1)}(Z_i,Z_j), \quad V_\nb^{(1)} :=\sum_{i<j} \phi_\nb^{(2)}(Z_i,Z_j), \quad W_\nb^{(1)} := \sum_{i<j} \phi_\nb^{(3)}(Z_i,Z_j), \notag
\end{align}
and define
\beq
U_{\nb}:= \frac{a_{0}}{l^2} M_\nb^{(1)} + \frac{a_{2}}{m^2} V_\nb^{(1)} +\frac{a_{1}}{2lm} W_\nb^{(1)} = \sum_{i<j}f_{\nb}(Z_i,Z_j). \notag
\eeq
With this definition of $U_\nb$, it follows from \eqref{jj:QnA} that decomposition \eqref{jj:dec} holds with remainder
\begin{eqnarray}
R_\nb   &=& R^{(1)}_\nb +  R^{(2)}_\nb + R^{(3)}_\nb + R^{(4)}_\nb + R^{(5)}_\nb + R^{(6)}_\nb, \notag
\end{eqnarray}
where
\begin{align*}
R_\nb^{(1)}                &:=  a_{0}\left ( \binom{n_1}{2}^{-1} - l^{-2} \binom{n_3}{2}^{-1} \right ) n\ep^{d/2} \bep^{-1}M_\nb^{(1)},  \\
R_\nb^{(2)}                &:= a_{2}\left ( \binom{n_2}{2}^{-1} - m^{-2} \binom{n_3}{2}^{-1} \right ) n\ep^{d/2}\bep^{-1} V_\nb^{(1)}, \\
R_\nb^{(3)}                &:= a_{1}\left ( (n_1 n_2)^{-1} - (2lm)^{-1} \binom{n_3}{2}^{-1} \right ) n\ep^{d/2}\bep^{-1} W_\nb^{(1)}, \\
R_\nb^{(4)}                &:= a_{0}\binom{n_1}{2}^{-1} n\ep^{d/2}\bep^{-1}M_\nb^{(2)},  \\
R_\nb^{(5)}                &:= a_{2}\binom{n_2}{2}^{-1} n\ep^{d/2}\bep^{-1}V_\nb^{(2)},  \\
R_\nb^{(6)}                &:=  a_{1}(n_1 n_2)^{-1} n\ep^{d/2} \bep^{-1}W_\nb^{(2)}.
\end{align*}
By the conventional theory of one-sample $U$-statistics (see, e.g., Ch.\ 1, Lee, 1990),
\beq\label{jj:rem2}
\Var(M_\nb^{(1)}) = \binom{n_3}{2}   \left (2(n_3-2)\xi_{1,\nb} + \xi_{2,\nb} \right),
\eeq
where
\beq
\xi_{1,\nb}:=\Cov(\phi_\nb^{(1)}(Z_1,Z_2),\phi_\nb^{(1)}(Z_1,Z_3)), \quad \xi_{2,\nb} := \Var(\phi_\nb^{(1)}(Z_1,Z_2)). \notag
\eeq
We have
\beq
\xi_{2,\nb} \leq \Ex(\phi_\nb^{(1)}(Z_1,Z_2)^2) \leq l^2 \Ex(\phi_\nb^{(1)}(Z_1,Z_2)) = l^4 P(d(X_1,X_2)<\ep), \notag
\eeq
where we use the fact that $\phi_\nb^{(1)}(Z_1,Z_2) \leq l^2$.
Hence $\xi_{1,\nb},\xi_{2,\nb}  =  \mbox{O}(\bep)$ and thus \eqref{jj:rem2} yields
\beq\label{jj:rem3}
\Var(M^{(1)}_\nb) = \mbox{O}(n_3^{3}\bep) \nbinf.
\eeq
Further, since
\beq\label{jj:coeff}
\binom{n_1}{2}^{-1} - l^{-2} \binom{n_3}{2}^{-1} = \binom{n_1}{2}^{-1} - l^{-2} \binom{n_1/l}{2}^{-1} \sim \frac{C}{n_1^3} \nbinf, \notag
\eeq
it follows from \eqref{jj:rem3} that
\beq\label{R1}
\Var(R_\nb^{(1)}) = \mbox{O}(n^{-1}) \to 0 \nbinf.
\eeq
Similarly, for $i = 2, 3$, we get
\beq
\Var(R_\nb^{(i)}) \to 0 \nbinf.
\eeq
Moreover, if a kernel is defined as
$$
\theta_\nb(z_i) := \sum_{1\leq j<k\leq l} I \left(d \left(x_{l(i-1)+j},x_{l(i-1)+k} \right)<\ep \right), \quad i = 1, \ldots, n_3,
$$
then
\beq
M_\nb^{(2)} = \sum_{i=1}^{n_3} \theta_\nb(Z_i) \notag
\eeq
and Lemma \ref{lemma} gives
\beq
\Var(M_\nb^{(2)}) = n_3 \Var(\theta_\nb(Z_1))\sim n_3 \binom{l}{2} \bep q_{2,0} = \mbox{O}(n_3\ep^d) \nbinf. \notag
\eeq
This implies
\beq
\Var(R_\nb^{(4)}) = \mbox{O}(n^{-1}) \to 0 \nbinf .
\eeq
In a similar way, for i=5, 6, we obtain
\beq\label{R56}
\Var(R_\nb^{(i)})  \to 0 \nbinf.
\eeq
Since $\Ex (R_\nb) = 0$, it follows from \eqref{R1}-\eqref{R56} that
\beq\label{jj:Rn}
R_\nb \prob 0 \nbinf.
\eeq
Next we prove asymptotic normality for $U_\nb$. Let
\beq\label{sig2:def}
\sigma^2_{\nb}  :=  \frac{n_3^2}{2} \Var(f_{\nb}(Z_1,Z_2)) + n_3^3 \Var( g_{\nb}(Z_1)).
\eeq
By applying Lemma \ref{lemma}, it is straightforward to show that,$\nbinf $,
\begin{align}\label{wn}
\sigma_\nb^2  \sim &  \frac{n_3^2 \bep}{2}   \Bigg(\frac{a_{0}^2}{l^2}q_{2,0}+ \frac{a_{2}^2}{m^2}q_{0,2} + \frac{a_{1}^2}{2lm}q_{1,1} \Bigg)  \\
& +  n_3^3 \bep^2\left (\frac{1}{l}\Var \left(a_{0}p_X(X)  +  \frac{a_{1}}{2}p_Y(X)\right)
 + \frac{1}{m} \Var \left(a_{2}p_Y(Y) + \frac{a_{1}}{2}p_X(Y) \right)\right). \notag
\end{align}
Since $n^2 \ep^d \to \infty$ leads to $n_3^2\bep \to \infty$, we get from \eqref{wn} that
$$
\sigma_\nb \to \infty \nbinf
$$
and thus
\beq\label{jj:c1}
\sup_{z_1,z_2} |f_{\nb}(z_1,z_2)| \leq  |a_{0}| + |a_{2}| + |a_{1}| = \mbox{o}(\sigma_{\nb}) \nbinf.
\eeq
Moreover, note that
\begin{align}\label{jj:c2}
\Ex (|f_{\nb}(z_1,Z_2)|)   \leq  \frac{|a_{0}|}{l} \sum_{i=1}^l p_{X,\ep}(x_i) +
\frac{|a_{2}|}{m} \sum_{i=1}^m p_{Y,\ep}(y_i) + \frac{|a_{1}|}{2l}\sum_{i=1}^l p_{Y,\ep}(x_i) + \frac{|a_{1}|}{2m} \sum_{i=1}^m p_{X,\ep}(y_i).
\end{align}
H\"{o}lder's inequality gives
\begin{align}\label{jj:unif}
p_{X,\ep}(x)  &= \int_{|y-x|<\ep} p_X(y)dy  \leq  \left(\int_{|y-x|<\ep} dy \right)^{1/2} \left(\int_{|y-x|<\ep} p_X(y)^2dy \right)^{1/2} \notag \\
 &= \bep^{1/2} \left(\int_{|y-x|<\ep} p_X(y)^2dy \right)^{1/2}, \notag
\end{align}
where the last integral tends to $0$ uniformly in $x$ as $\ep \to 0$.
The corresponding results can be shown for the other terms in \eqref{jj:c2}.
Hence, we obtain
\beq\label{jj:c3}
\sup_{z_1}\Ex( |f_{\nb}(z_1,Z_2)|)  = \mbox{o}(\bep^{1/2}) = \mbox{o}(\sigma_{\nb}/n_3) \nbinf.
\eeq
Now it follows from  \eqref{jj:c1} and \eqref{jj:c3} that the conditions of Theorem 2.1 in Jammalamadaka and Janson (1986) are fulfilled.
Consequently, it holds that
\beq\label{jj:normal}
\frac{U_{\nb} - \Ex (U_{\nb}) }{\sigma_{\nb}} = \frac{1}{\sigma_{\nb}} \left (\sum_{i < j}f_{\nb}(Z_i,Z_j) - \binom{n_3}{2} \mu_\nb \right )\weak N(0,1) \nbinf.
\eeq
Furthermore, since $n\ep^d \to \beta, 0\leq \beta < \infty$, we get from \eqref{wn} and Lemma \ref{lemma} that
\beq\label{jj:varas}
n^2 \ep^d \binom{n_3}{2}^{-2} \bep^{-2} \sigma_\nb^2 \to  \eta +  \beta \zeta \nbinf.
\eeq
Finally, assertion \eqref{nep^d/2} follows from \eqref{jj:dec}, \eqref{jj:Rn}, \eqref{jj:normal}, \eqref{jj:varas}, and the Slutsky theorem. This completes the proof. \qed
\\\\\\
{\it Proof of Theorem \ref{th:main3}.} \textit{(i)}
The proof is similar to that of Theorem 7 in Leonenko and Seleznjev (2010).
Let $V := (V_1,\ldots,V_d)'$ be an auxiliary random vector uniformly distributed in the unit ball $B_1(0)$.
By definition, we have $\qt_{1,1,\ep} = \bep^{-1} \Ex(p_{X,\ep}(Y)) = \Ex(p_X(Y-\ep V))$ and thus
\begin{align}
\qt_{1,1,\ep}-q_{1,1}  =&\   \int_{R^d} \Ex ( p_X(y-\ep V))p_Y(y)dy - \int_{R^d} p_X(y)p_Y(y)dy \notag \\
=& \   \Ex \left( \int_{R^d} ( p_X(y-\ep V) - p_X(y))p_Y(y)dy \right) \notag \\
 =&\   \Ex \left( \int_{R^d} ( p_X(y-\ep V) - p_X(y))(p_Y(y)- p_Y(y-\ep V))dy \right) \notag \\
&+   \Ex \left (  \int_{R^d} ( p_X(y-\ep V) - p_X(y))p_Y(y-\ep V)dy \right). \notag
\end{align}
For the last term, using the change of variables $z=y -\ep V$ and symmetry $V \stackrel{\mathrm{D}}{=} -V$, we obtain
\begin{align}
 &\Ex \left( \int_{R^d} (p_X(y-\ep V) -  p_X(y))p_Y(y-\ep V)dy \right) \notag \\
 & = \Ex \left( \int_{R^d} (p_X(z) - p_X(z+\ep V))p_Y(z)dz \right)  \notag \\
&= \Ex \left( \int_{R^d} (p_X(z) - p_X(z-\ep V))p_Y(z)dz \right) = -(\qt_{1,1,\ep}-q_{1,1}). \notag
\end{align}
From the above,
\beq
2(\qt_{1,1,\ep}-q_{1,1}) = \Ex \left ( \int_{R^d} ( p_X(y-\ep V) - p_X(y))(p_Y(y)- p_Y(y-\ep V))dy \right)\notag
\eeq
and hence H\"{o}lder's inequality and the density smoothness condition imply
\begin{align*}
| \qt_{1,1,\ep}-q_{1,1}| & \leq  \frac{1}{2} \Ex  \left(\int_{R^d}(p_X(y-\ep V) - p_X(y))^2dy  \right )^{1/2}
\left ( \int_{R^d}(p_Y(y)- p_Y(y-\ep V))^2dy \right)^{1/2}   \\
 & \leq \frac{1}{2}K^2 \Ex(|V|^{2\alpha}) \ep^{2\alpha} \leq \frac{1}{2} K^2 \ep^{2\alpha}. \notag
\end{align*}
Similar inequalities can be obtained for $ \qt_{2,0,\ep}$ and  $ \qt_{0,2,\ep}$.
It follows that
\beq
|\lsqep-\lsq| \leq C \ep^{2\alpha}, \quad C := \frac{1}{2} K^2(|a_{0}|+ |a_{1}| + |a_{2}|), \notag
\eeq
and so the assertion is proved.
\\\\
\textit{(ii)} First note that the conditions $\ep \sim cn^{-1/(2\alpha+d/2)},c>0$, and $0 < \alpha \leq d/4$ yield
\beq\label{n2ed}
n^2\ep^d \sim c^d n^{\frac{4\alpha}{2\alpha + d/2}} \leq c^d n \nbinf.
\eeq
For $\kb \in \{(2,0), (1,1),(0,2) \}$, by combining \eqref{eq:var}, \eqref{eq:bound}, \eqref{eq:varlim}, and \eqref{n2ed}, we get that
\begin{equation*}
\Var(\Qt_{\kb,\nb}) = \Or\left( \frac{1}{n^2\ep^d} \right) = \Or\left( \frac{1}{n^{4\alpha/(2\alpha + d/2)}} \right) \nbinf
\end{equation*}
and consequently
\begin{align}\label{eq:vn}
\Var(\Lsq) = \Or\left( \frac{1}{n^{4\alpha/(2\alpha + d/2)}} \right) \nbinf.
\end{align}
Moreover, from \textit{(i)} we have $(\lsqep-\lsq)^2 \lesssim C_1 n^{-\frac{4\alpha}{2\alpha + d/2}},C_1>0$, which together with \eqref{eq:vn} gives
\beq
\Var(\Lsq) +(\lsqep-\lsq)^2 = \Or \left ( \frac{1}{n^{4\alpha/(2\alpha + d/2)}}  \right) \nbinf. \notag
\eeq
Using this, for some $C_2>0$, any $A>0$, and large enough $n_1,n_2$, we obtain
\begin{align}
P\left(|\Lsq-\lsq | >An^{-\frac{2\alpha}{2\alpha + d/2}}  \right )& \leq n^{\frac{4\alpha}{2\alpha + d/2}} \frac{\Var(\Lsq) + (\lsqep-\lsq)^2}{A^2} \leq \frac{C_2}{A^2} \notag
\end{align}
and the assertion follows.\\\\
\textit{(iii)} The argument is similar to that of {\it (ii)} and therefore is left out. This completes the proof. \qed
\\\\\\
{\it Proof of Theorem \ref{th:main4}.} \textit{(i)}
We have
\beq\label{eq:nepd}
\sqrt{n}(\Lsq - \lsq) = \sqrt{n}(\Lsq-\lsqep) + \sqrt{n}(\lsqep-\lsq),
\eeq
where, by Theorem \ref{th:main3}{\it (i)} in combination with the conditions
$\ep \sim L(n)n^{-1/d}$ and $\alpha > d/4$,
\beq
|\sqrt{n}(\lsqep-\lsq)| \leq  C n^{1/2} \ep^{2\alpha} \sim
C L(n)^{2\alpha} n^{1/2-2\alpha/d} \to 0 \nbinf. \notag
\eeq
The assertion thus follows from Theorem \ref{th:main2}\textit{(i)}, \eqref{eq:nepd}, and the Slutsky theorem.\\\\
\noindent{\it (ii)} Consider the decomposition corresponding to \eqref{eq:nepd}:
\begin{equation}  \label{eq:Norm2:dec}
n^{\gamma/(1+\gamma)} c^{d/2}(\Lsq- \lsq) = n^{\gamma/(1+\gamma)} c^{d/2}(\Lsq -
\lsqep) + n^{\gamma/(1+\gamma)} c^{d/2}(\lsqep-\lsq).
\end{equation}
Since $\ep \sim cn^{-2/((1+\gamma)d)},c>0$, for some $0<\gamma<1$, we get $n\ep^d \to 0$, $n^2 \ep^d \to \infty$, and $n^{\gamma/(1+\gamma)} c^{d/2} \sim n\ep^{d/2}$,
so Theorem \ref{th:main2}\textit{(ii)} implies the asymptotic normality
\begin{equation}\label{eq:Norm2:norm}
n^{\gamma/(1+\gamma)} c^{d/2}(\Lsq-\lsqep) \overset{%
\mathrm{D}}{\to} N(0, \eta) \nbinf.
\end{equation}
Further, Theorem \ref{th:main3}{\it (i)} and the assumptions $\ep \sim cn^{-2/((1+\gamma)d)},c>0$, and $\alpha > (d/4)\gamma$ lead to
\begin{equation}  \label{eq:Norm2:bias}
|n^{\gamma/(1+\gamma)} c^{d/2}(\lsqep-\lsq)| \leq
c^{d/2} C n^{\gamma/(1+\gamma)}\ep^{2\alpha} \sim c^{d/2+2\alpha}Cn^{\gamma/(1+\gamma) - 4\alpha/((1+\gamma)d)} \to 0 \nbinf,
\end{equation}
where the last limit holds because $\gamma/(1+\gamma) - 4\alpha/((1+\gamma)d) < 0$.
Now the statement follows from \eqref{eq:Norm2:dec}, \eqref{eq:Norm2:norm}, \eqref{eq:Norm2:bias},
and the Slutsky theorem. \\\\
\noindent{\it (iii)} 
From Theorem \ref{th:main3}\textit{(i)} and the condition $\ep \sim L(n)^{2/d}n^{-2/d}$, we obtain
\beq\label{L(n)}
|L(n)(\lsqep-\lsq)| \leq C L(n)\ep^{2\alpha} \sim  C L(n)^{1+4\alpha/d}n^{-4\alpha/d} \to 0 \nbinf.
\eeq
Note also that $\ep \sim L(n)^{2/d}n^{-2/d}$ gives $n\ep^d \to 0$ and $L(n) \sim  n\ep^{d/2}$.
Therefore, similarly as above, the assertion is implied by the decomposition corresponding to \eqref{eq:Norm2:dec}, \eqref{L(n)},
Theorem \ref{th:main2}\textit{(ii)}, and the Slutsky theorem.
This completes the proof.
\qed\\\\\\
{\it Proof of Proposition \ref{test}.} \textit{(i)} When $n^2\ep^d \to \infty$ and $n\ep^{d} \to \beta, 0 \leq \beta < \infty$,
Theorem~\ref{th:main2} can be applied with $\Lsq = \hat{D}_{2,\nb}$.
Indeed, under $H_0$ we have $\lsqep=\Ex(\hat{D}_{2,\nb}) = 0$ and $\zeta = 0$, so Theorem \ref{th:main2} yields
\beq\label{Tn}
n \ep^{d/2}\hat{D}_{2,\nb} \weak N(0, \eta) \nbinf \notag
\eeq
in this case.
Hence, we need to show that, for $\hat{D}_{2,\nb}$ under $H_0$, the proof of Theorem \ref{th:main2}
can be modified so that the assumption $n\ep^d \to \beta, 0 \leq \beta < \infty$, is unnecessary.
In fact, this assumption is only needed for convergence property \eqref{jj:varas} of $\sigma_\nb^2$.
Under $H_0$, we obtain from definition \eqref{jj:fngn} that $g_\nb(z)=0$ and thus $\Var(g_\nb(Z_1)) = 0$.
Therefore, definition \eqref{sig2:def} of $\sigma_\nb^2$ implies
\beq
\sigma_\nb^2 \sim \; \frac{n_3^2 \bep}{2}   \Bigg(\frac{a_{0}^2}{l^2}q_{2,0}+ \frac{a_{2}^2}{m^2}q_{0,2} + \frac{a_{1}^2}{2lm}q_{1,1} \Bigg) \nbinf \notag
\eeq
and hence \eqref{jj:varas} can be written
\beq
n^2 \ep^d \binom{n_3}{2}^{-2} \bep^{-2} \sigma_\nb^2 \to \eta \nbinf, \notag
\eeq
which does not require convergence of $n\ep^d$.
The assertion follows.
\\\\
\noindent
\textit{(ii)} Under $H_1$, we get from Theorem \ref{th:main1} and the Slutsky theorem that
$\hat{D}_{2,\nb}/v_\nb \stackrel{P}{\to} D_2/\sqrt{\eta} > 0$.
Consequently, since $c_\nb = \oor(n\ep^{d/2})$, we see that
\beq
P(T_\nb > c_\nb) = P(\hat{D}_{2,\nb}/v_\nb > c_\nb /(n\ep^{d/2})) \to 1 \nbinf. \notag
\eeq
This completes the proof.
\qed\\\\\\
{\it Proof of Proposition \ref{log}.}
The assertion follows straightforwardly from Theorem \ref{th:main4} in a similar way as in Leonenko and Seleznjev (2010).
\qed
\bigskip

\subsection*{Acknowledgment}

{}
The second author is partly supported by the Swedish Research Council grant
2009-4489. We  would like to thank  professor
N.\ Leonenko for valuable comments and  discussions.\\
\medskip

{\noindent {\large \textbf{References}}}

\begin{reflist}

 {\small Ahmad, I.A., Cerrito, P.B., 1993, Goodness of fit tests based on the $L_2$-norm of multivariate probability density functions,
 {J.\ Nonparametr.\ Stat.} 2, 169--181.}

{\small Basseville, M., 2010, Divergence measures for statistical data processing,
Technical Report 1961, IRISA.}

 {\small Basu, A., Harris, I.R., Hjort, N.L., Jones, M.C., 1998, Robust and efficient estimation by minimising a density power divergence,
 {Biometrika} 85, 549--559.}

 {\small Bickel, P.J. and Ritov, Y., 1988, Estimating integrated squared density derivatives: sharp best order of convergence estimates,
{Sankhy{\=a}: The Indian Journal of Statistics} Series A, 381--393.}

{\small Bogachev, V.I., 2007, Measure Theory, vol. I, Springer-Verlag, Berlin.}

 {\small Broniatowski, M., Toma, A., Vajda, I., 2012, Decomposable pseudodistances and applications in statistical estimation,
{J.\ Statist.\ Plann.\ Inference}, 142, 2574--2585.}

{\small Escolano, F., Suau, P., Bonev, B., 2009, Information Theory in Computer Vision and Pattern Recognition, Springer, New York.}

 {\small Gin\'{e}, E., Nickl, R., 2008, A simple adaptive estimator for the integrated square of a density,
 {Bernoulli} 14, 47--61.}

{\small Jammalamadaka, S.R., Janson, S., 1986,
Limit theorems for a triangular scheme of $U$-statistics with applications to inter-point distances,
{Ann.\  Probab.} 14, 1347--1358.}

{\small Kallenberg, O., 1997, Foundations of Modern Probability, Springer-Verlag, New York.}

{\small K\"{a}llberg, D., Leonenko, N., Seleznjev, O., 2012, Statistical inference for R\'{e}nyi entropy functionals,
{Lecture Notes in Comput.\ Sci.}\ 7260, 36--51.}

{\small Laurent, B., 1996, Efficient estimation of integral functionals of a density,
{Ann.\ Statist.} 24, 659--681.}

{\small Lee, A.J., 1990, $U$-Statistics: Theory and Practice, Marcel Dekker, New York.}

{\small Leonenko,  N., Pronzato, L., Savani, V., 2008, A class of R\'{e}nyi information estimators for multidimensional densities.
{Ann.\  Statist.}\ 36, 2153--2182. Corrections, 2010, {Ann.\ Statist.}\ 38, 3837--3838.}

{\small  Leonenko, N., Seleznjev, O., 2010, Statistical inference for the $\ep$-entropy and the quadratic R\'{e}nyi entropy.
{J.\  Multivariate Anal.}\ 101, 1981--1994.}

 {\small Li, Q., 1996, Nonparametric testing of closeness between two unknown distribution functions,
 {Econometric Rev.}\ 15, 261--274.}

 {\small Neemuchwala, H., Hero, A., Carson, P., 2005, Image matching using alpha-entropy measures and entropic graphs,
 {Signal Processing} 85, 277--296.}

{\small Pardo, L., 2006, Statistical Inference Based on Divergence Measures,
Chapman \& Hall, Boca Raton.}

{\small Penrose, M., 1995, Generalized two-sample $U$-statistics and a two-species reaction-diffusion model,
{Stochastic Process.\ Appl.}\ 55, 57--64.}

{\small Principe, J.C., 2010, Information Theoretic Learning, Springer, New York.}

{\small R\'{e}nyi, A., 1970, Probability Theory, North-Holland, Amsterdam.}

{\small Seleznjev, O., Thalheim, B., 2003, Average case analysis in database problems,
{Methodol.\ Comput.\ Appl.\ Prob.}\ 5, 395--418. }

 {\small Seleznjev, O., Thalheim, B., 2010, Random databases with approximate record matching,
 {Methodol.\ Comput.\ Appl.\ Prob.}\ 12, 63--89. }

{\small Thalheim, B., 2000, Entity-Relationship Modeling. Foundations of Database Technology, Springer-Verlag, Berlin.}

{\small Ullah, A., 1996, Entropy, divergence and distance measures with econometric applications.
{J. Statist. Plann. Inference} 49, 137--162.}

 {\small Weber, N.C., 1983, Central limit theorems for a class of symmetric statistics,
 {Math.\ Proc.\ Cambridge Philos.\ Soc.}\ 94, 307--313. }

 \end{reflist}

\end{document}